\documentclass{article}
\usepackage{amssymb}
\usepackage{amsmath}
\usepackage{epsfig}
\usepackage{color}
\usepackage{epic,eepic}

\author{Sander Dahmen}
\title{Counting Integral Lam\'e Equations by Means of Dessins d'Enfants}
\date{November 27, 2003}
\newtheorem{theorem}{Theorem}

\newtheorem{corollary}[theorem]{Corollary}

\newenvironment{proof}{\textbf{Proof}}{\hfill $\square$ \\}
\newcommand{\Set}[2]{\ensuremath{\{#1 \: | \; #2 \}}}
\newcommand{\N}{\ensuremath{\mathbb{N}}}
\newcommand{\Z}{\ensuremath{\mathbb{Z}}}
\newcommand{\Q}{\ensuremath{\mathbb{Q}}}

\newcommand{\C}{\ensuremath{\mathbb{C}}}

\begin{document}

\maketitle

\begin{abstract}
We obtain an explicit formula for the number of Lam\'e equations
(modulo scalar equivalence) with index $n$ and projective
monodromy group of order $2N$, for given $n \in \Z$ and $N \in
\N$. This is done by performing the combinatorics of the `dessins
d'enfants' associated to the Belyi covers which transform
hypergeometric equations into Lam\'e equations by pull-back.
\end{abstract}

\section{Introduction}

The integral Lam\'e equation with parameters $B,g_2,g_3 \in \C$
and $n \in \Z$ is the second order differential equation on the
$\mathbb{P}^1$ given by
\begin{equation}\label{eqn Lame}
p(z)\frac{d^2y}{dz^2} + \frac{1}{2} p'(z) \frac{dy}{dz} - (n(n+1)z
+ B)y = 0,
\end{equation}
where $p(z):= 4z^3 -g_2z -g_3$ has nonzero discriminant, i.e.
$g_2^3-27 g_3^2\not=0$. The parameter $n$ is called the
\emph{index}. We are interested in Lam\'e equations with a basis
of solutions which are algebraic (over $\C(z)$). Such solutions
can in fact occur for non integer $n \in \Q$, but in this article
we restrict to $n \in \Z$, hence the word integral in `integral
Lam\'e equation'.

It is known that if there exists a basis of algebraic solutions,
then the projective monodromy group is dihedral, see for example
\cite[Cor. 3.4]{Beukers vd Waal} or \cite[Thm. 4.5.5]{vd Waal}. If
in equation (\ref{eqn Lame}) we replace $z$ by $\lambda z$ (for a
nonzero $\lambda \in \C$), we get a new Lam\'e equation with the
same index $n$, but with the parameters $(B,g_2,g_3)$ replaced by
$(B/\lambda,g_2/\lambda^2,g_3/\lambda^3)$. These substitutions
induce a natural equivalence relation on the space of all Lam\'e
equations. Two Lam\'e equations which are equivalent w.r.t. this
equivalence relation are called \emph{scalar equivalent}. Now it
is also known that for given $n \in \Z$ and $N \in \N$ there are
only finitely many Lam\'e equations modular scalar equivalence
with index $n$ and projective monodromy group dihedral of order
$2N$. This can be proven in different ways, see for example
\cite[Thm. 4.6]{Beukers vd Waal}, \cite[Thms. 5.4.4, 6.7.9]{vd
Waal} or \cite[Thm 4.1]{Litcanu 1}. Throughout this article the
number of Lam\'e equations modular scalar equivalence with index
$n$ and projective monodromy group dihedral of order $2N$ will be
denoted by $L(n,N)$.

In \cite{Litcanu 1} it is described how the problem of calculating
$L(n,N)$ can be translated in counting the number of dessins
d'enfants compatible with prescribed ramification data. In the
same article the combinatorics are performed for $n=1$, obtaining
a result earlier obtained in \cite{Chiarellotto}. In \cite{Litcanu
2} an attempt was made to perform the combinatorics for $n=2$. In
this article we perform the combinatorics for general $n \in \N$.

\section{The Combinatorics}

Following \cite{Litcanu 1}, we write the Lam\'e equation as in
equation \ref{eqn Lame}, but now with $p(z):=4z(z-1)(z-\lambda)$
and $\lambda \in \C-\{0,1\}$. According to \cite{Chiarellotto} the
functions (Belyi covers) $F:\mathbb{P}^1 \to \mathbb{P}^1$ which
transform, by pull-back, hypergeometric equations into Lam\'e
equations with index $n$ and projective monodromy group dihedral
of order $2N$, have the following ramification data.

\begin{table}[h!]
\begin{tabular}{c|c|c|c|c|c}
& 0 & 1 & $\lambda$ & $\infty$  & \\
\hline
1 & & & & & $+nN/2$ points with multiplicity 2 \\
\hline
$\infty$ & & & & & $+n$ points with multiplicity N \\
\hline
0 & 1 & 1 & 1 & 2n+1 & $+(nN - 2n -4)/2$ points with multiplicity 2 \\
\end{tabular}
\caption{Case Ia}\label{tabel Ia}
\end{table}

\begin{table}[h!]
\begin{tabular}{c|c|c|c|c|c}
& 0 & 1 & $\lambda$ & $\infty$  & \\
\hline
1 & & & 1 & & $+(nN-1)/2$ points with multiplicity 2 \\
\hline
$\infty$ & & & & & $+n$ points with multiplicity N \\
\hline
0 & 1 & 1 & & 2n+1 & $+(nN - 2n -3)/2$ points with multiplicity 2 \\
\end{tabular}
\caption{Case Ib}\label{tabel Ib}
\end{table}

\begin{table}[h!]
\begin{tabular}{c|c|c|c|c|c}
& 0 & 1 & $\lambda$ & $\infty$  & \\
\hline
1 & & 1 & 1 & & $+(nN-2)/2$ points with multiplicity 2 \\
\hline
$\infty$ & & & & & $+n$ points with multiplicity N \\
\hline
0 & 1 & & & 2n+1 & $+(nN - 2n -2)/2$ points with multiplicity 2 \\
\end{tabular}
\caption{Case Ic}\label{tabel Ic}
\end{table}

\begin{table}[h!]
\begin{tabular}{c|c|c|c|c|c}
& 0 & 1 & $\lambda$ & $\infty$  & \\
\hline
1 & 1 & 1 & 1 & & $+(nN-3)/2$ points with multiplicity 2 \\
\hline
$\infty$ & & & & & $+n$ points with multiplicity N \\
\hline
0 & & & & 2n+1 & $+(nN - 2n -1)/2$ points with multiplicity 2 \\
\end{tabular}
\caption{Case Id}\label{tabel Id}
\end{table}

\begin{table}[h!]
\begin{tabular}{c|c|c|c|c|c}
& 0 & 1 & $\lambda$ & $\infty$  & \\
\hline
1 & & & & & $+nN/2$ points with multiplicity 2 \\
\hline
$\infty$ & & N/2 & N/2 & & $+n-1$ points with multiplicity N \\
\hline
0 & 1 & & & 2n+1 & $+(nN - 2n -2)/2$ points with multiplicity 2 \\
\end{tabular}
\caption{Case II}\label{tabel II}
\end{table}

The tables have the following meaning. All the possible branching
points, i.e. $0,1$ and $\infty$, are contained in the first
column. In the three corresponding rows, the inverse images (which
are given in the first row) together with the multiplicities can
be read of in the obvious way. For more information we refer to
\cite{Litcanu 1} or \cite{Litcanu 2}.

\begin{theorem}\label{thm number of dessins}
Let $n,N \in \N$, then the number of dessins d'enfants compatible
with the tables above equals $n(n+1)(N-1)(N-2)/12 + 2/3
\varepsilon(n,N)$, where
\[ \varepsilon(n,N):=\left\{
\begin{array}{ll}
1 & \mbox{if $3|N$ and $n \equiv 1 \pmod{3}$}; \\
0 & \mbox{otherwise}.
\end{array}
\right. \]
\end{theorem}

\begin{proof}
Let $n,N \in \N$ and suppose $n > 1$. First consider case I, it
will not be necessary to consider the four cases Ia, Ib, Ic, Id
separately. In all four cases there are $n$ points above $\infty$,
all with multiplicity $N$. So the associated dessins consist of
$n$ cells all with valency $N$. Furthermore, there are $n-1$
cycles, which have  beginning and end in the point $p:=\infty$ and
3 lines, emanating from $p$. There are no further intersections.
In this article, whenever we draw dessins d'enfants, we will not
draw the vertices.

\begin{center}
\begin{picture}(0,0)%
\includegraphics{linecycle.pstex}%
\end{picture}%
\setlength{\unitlength}{3947sp}%
\begingroup\makeatletter\ifx\SetFigFont\undefined%
\gdef\SetFigFont#1#2#3#4#5{%
  \reset@font\fontsize{#1}{#2pt}%
  \fontfamily{#3}\fontseries{#4}\fontshape{#5}%
  \selectfont}%
\fi\endgroup%
\begin{picture}(2450,1180)(5464,-3940)
\end{picture}

\begin{figure}[h!]
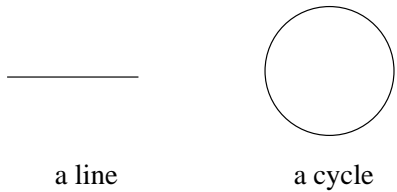

\caption{a line and a cycle, not taken into account the valencies}
\label{fig linecycle}
\end{figure}
\end{center}

We define the valency of a line as the number of edges of that
line and the valency of a cycle as half the number of edges of
that cycle. With these definitions we define the valency of a cell
as the sum of the valencies of the lines it contains and the
cycles that bound it. We note that the valency of a cell equals
the multiplicity of the inverse image of $\infty$ lying in the
cell. We see that the boundary of every cell contains exactly one
ore two cycles (since we assumed $n > 1$). A cell with exactly one
cycle in its boundary will be called a \emph{simple cell}. Since
all the cells have the same valency and there are only three
lines, there can be no more than three simple cells. There is of
course a minimum of two simple cells. These two cases, that of two
simple cells (case A) and that of three simple cells (case B),
will be considered seperately.

We first consider case A. If we do not take into account the 3
lines and the valencies of the cycles, it is easy to see that
there is only one possibility for the dessin. It has the following
shape:

\begin{center}
\begin{picture}(0,0)%
\includegraphics{caseA.pstex}%
\end{picture}%
\setlength{\unitlength}{3947sp}%
\begingroup\makeatletter\ifx\SetFigFont\undefined%
\gdef\SetFigFont#1#2#3#4#5{%
  \reset@font\fontsize{#1}{#2pt}%
  \fontfamily{#3}\fontseries{#4}\fontshape{#5}%
  \selectfont}%
\fi\endgroup%
\begin{picture}(1516,1663)(5243,-3668)
\put(6001,-2161){\makebox(0,0)[lb]{\smash{\SetFigFont{12}{14.4}{\rmdefault}{\mddefault}{\updefault}{\color[rgb]{0,0,0}$n-1$ cycles}%
}}}
\end{picture}

\begin{figure}[h!]
\caption{case A}
\label{fig caseA}
\end{figure}
\end{center}

This can also be drawn in a more symmetric way as follows:

\begin{center}
\input{symcaseA.pstex_t}
\begin{figure}[h!]
\caption{case A}
\label{fig symcaseA}
\end{figure}
\end{center}

We will now take into account the three lines (but still do not
take into account the valencies). Each of the two simple cells
must contain at least one line (since every cell has the same
valency). If each of the three lines are contained in a simple
cell, there is (because of rotational symmetry) only one
possibility. We call this case AI. The dessin has the following
shape:

\begin{center}
\input{caseAI.pstex_t}
\begin{figure}[h!]
\caption{case AI}
\label{fig caseAI}
\end{figure}
\end{center}

If there is a line which is not contained in a simple cell, there
are $n-2$ possible cells that can contain it. For every of these
$n-2$ cells the line can lie on two different sides, but if we
take into account the rotational symmetry, we see that the number
of possibilities is reduced by a factor two. So in this case we
arrive at exactly $n-2$ different kinds of dessins. We call this
case AII. The dessins are of the following shape:

\begin{center}
\input{caseAII.pstex_t}
\begin{figure}[h!]
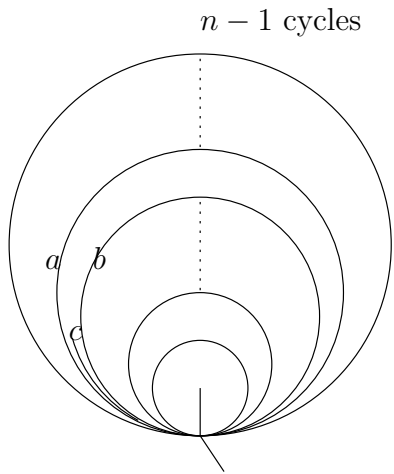

\caption{case AII}
\label{fig caseAII}
\end{figure}
\end{center}

where the line which is not contained in a simple cell can always
be drawn on the left side.

We will now take into account the valencies of the cycles and
lines. Let $a,b,c \in \N$ denote the valencies as given in figures
\ref{fig caseAI} and \ref{fig caseAII}. In both cases we have
$a+b+c=N$ (and $a,b,c \geq 1$). Furthermore, it is easy to see (by
induction) that for every such triple $(a,b,c)$ there is exactly
one possible dessin (for a fixed shape of the dessin). Now there
are $(N-1)(N-2)/2$ possible triples $(a,b,c)$. We conclude that
there are exactly $(N-1)(N-2)/2$ possible dessins of type AI and
there are exactly $(n-2)(N-1)(N-2)/2$ possible dessins of type
AII. So there are exactly $(n-1)(N-1)(N-2)/2$ possible dessins of
type A.

Now consider case B. If we do not take into account the valencies
of the cells and lines, one easily obtains (by induction) that the
dessins are of the following shape:

\begin{center}
\input{caseB.pstex_t}
\begin{figure}[h!]
\caption{case B}
\label{fig caseB}
\end{figure}
\end{center}

We first consider these dessins without taking into account
rotational symmetry. Since $x+y+z=n-1$ (and $x,y,z \geq 1$) there
are $(n-2)(n-3)/2$ possibilities for the triple $(x,y,z)$. Note
that for $n=2,3$ we also arrive at the correct answer, namely
zero. We will now take into account the valencies of the cycles
and lines. Let $a,b,c \in \N$ denote the valencies as given in
figure \ref{fig caseB}. We have $a+b+c=N$ (and $a,b,c \geq 1$).
Furthermore, it is easy to see (by induction) that for every such
triple $(a,b,c)$ there is exactly one possible dessin (for fixed
$x,y,z$ and not taken into account the rotational symmetry yet).
Now there are $(N-1)(N-2)/2$ possible triples $(a,b,c)$. We
conclude that there are $(n-2)(n-3)(N-1)(N-2)/4$ possible dessins
of type B, not taken into account the rotational symmetry.

If we now do take into account the rotational symmetry, we see
that we have counted every dessin three times, except when the
dessin has a three fold rotational symmetry. There is exactly one
dessin with a three fold rotational symmetry if and only if $N$
and the number of cycles are divisible by 3, i.e. $3|N$ and $n
\equiv 1 \pmod{3}$, otherwise there are no such dessins. We
conclude that there are exactly
\[\frac{1}{3}\left(\frac{(n-2)(n-3)(N-1)(N-2)}{4} +
2\varepsilon(n,N) \right)=
\frac{(n-2)(n-3)(N-1)(N-2)}{12}+\frac{2}{3} \varepsilon(n,N)\]
possible dessins of type B.

So the total number of dessins of type I equals
$(n-1)(N-1)(N-2)/2+(n-2)(n-3)(N-1)(N-2)/12 + 2/3
\varepsilon(n,N)=n(n+1)(N-1)(N-2)/12+2/3 \varepsilon(n,N)$.

It remains to be proven that there are no dessins of type II. In
case II there are $n+1$ points above $\infty$, $n-1$ with
multiplicity $N$ and two with multiplicity $N/2$. So the
associated dessins consist of $n+1$ cells, $n-1$ with multiplicity
$N$ and two with multiplicity $N/2$. Furthermore, there are $n$
cycles and 1 line, which emanate from $p=\infty$. There are no
further intersections. From the fact that there are two cells with
valency strictly smaller than $N$, all other cells have valency
$N$ and there is one line, we conclude that (again) there are
either 2 simple cells (case C) or 3 simple cells (case D).

First consider the case of two simple cells. We distinguish the
following two cases. Case CI: the line is contained in a simple
cell. Case CII: the line is not contained in a simple cell.

\begin{center}
\input{caseC.pstex_t}
\begin{figure}[h!]
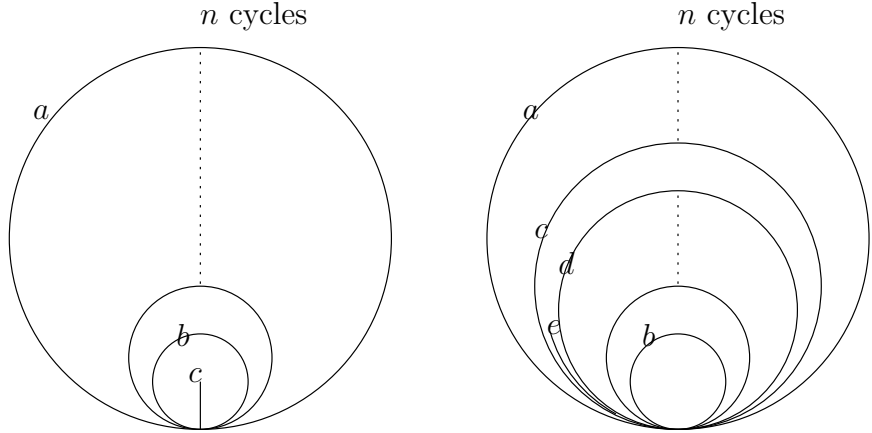

\caption{left: case CI, right: case CII }
\label{fig caseC}
\end{figure}
\end{center}

Let $a,b,c \in \N$ denote the valencies as given in the left part
of figure \ref{fig caseC}. We must have $a=N/2$ and by induction
we obtain that all the non simple cells have valency $N$ and that
$b=N/2$. So $b+c
>N/2$, but this means that there is only one cell with valency
$N/2$. We conclude that there are no dessins of type CI.

Let $a,b,c,d,e \in \N$ denote the valencies as given in the right
part of figure \ref{fig caseC}. We must have $a=b=N/2$. Since all
non simple cells must have valency $N$, we obtain by induction
that $c=d=N/2$. So $c+d+e
>N$, but this means that there is a cell with valency larger that
$N$. We conclude that there are no dessins of type CII.

Now consider the case of three simple cells. The line must lie in
a simple cell and the dessin has the following shape:

\begin{center}
\input{caseD.pstex_t}
\begin{figure}[h!]
\caption{case D}
\label{fig caseD}
\end{figure}
\end{center}

Let $a,b,c,d,e \in \N$ denote the valencies as given in figure
\ref{fig caseD}. We must have $a=b=N/2$. Since all non simple
cells must have valency $N$ we obtain by induction that $c=d=N/2$.
So $c+d+e>N$, but this means that there is a cell with valency
larger that $N$. We conclude that there are also no dessins of
type D. So there are no dessins of type II. We note that for $n=1$
the proof that there are no dessins of type II also applies (but
only case CI has to be considered and there is only one cycle).

We have proven our theorem for $n>1$. For $n=1$ the combinatorics
was done in \cite{Litcanu 1} (and for $n=1$ the result was already
obtained in \cite{Chiarellotto} by other means), but the
combinatorics can be simplified significantly as follows.

As noted before there are no dessins corresponding to type II. So
let us consider case I. There is only one point above $\infty$,
which has multiplicity $N$. So the associated dessins consist of
one cell with valency $N$. Furthermore there are no cycles and
three lines, who come together in the point $p=\infty$. The dessin
has the following shape:

\begin{center}
\begin{picture}(0,0)%
\includegraphics{case_n1.pstex}%
\end{picture}%
\setlength{\unitlength}{3947sp}%
\begingroup\makeatletter\ifx\SetFigFont\undefined%
\gdef\SetFigFont#1#2#3#4#5{%
  \reset@font\fontsize{#1}{#2pt}%
  \fontfamily{#3}\fontseries{#4}\fontshape{#5}%
  \selectfont}%
\fi\endgroup%
\begin{picture}(2424,1824)(4789,-4573)
\end{picture}

\begin{figure}[h!]
\caption{$n=1$}
\label{fig case_n1}
\end{figure}
\end{center}

Let $a,b,c$ denote the valencies of the three lines. If we do not
take into account the rotational symmetry, then the number of
dessins equals the number of triples $(a,b,c)$ with $a+b+c=N$ (and
$a,b,c \geq 1$). There are $(N-1)(N-2)/2$ such triples.

If we now do take into account the rotational symmetry, we see
that we have counted every dessin three times, except when the
dessin has a three fold rotational symmetry. There is exactly one
dessin with a three fold rotational symmetry if and only if $3|N$,
otherwise there are no such dessins. We conclude that there are
exactly
\[\frac{1}{3}\left(\frac{(N-1)(N-2)}{2} +
2\varepsilon(1,N) \right)= \frac{(N-1)(N-2)}{6}+\frac{2}{3}
\varepsilon(1,N)\] possible dessins when $n=1$. This finishes the
proof.
\end{proof}

Now theorem \ref{thm number of dessins} combined with proposition
3.1 in \cite{Litcanu 1} gives rise to the following.

\begin{theorem}
Let $n \in \Z$ and $N \in \N$. Then
\begin{equation}\label{eqn sum L(n,N)}
\sum_{d|N}L(n,d)= \frac{n(n+1)}{12}(N-1)(N-2) + \frac{2}{3}
\varepsilon(n,N),
\end{equation}
where $\varepsilon(n,N)$ is as in theorem \ref{thm number of
dessins}.
\end{theorem}

Let $N \in \N$. We denote by $\phi(N)$ Euler's totient function,
i.e.
\[\phi(N):=|\Set{k \in \Z}{0 \leq k < N,\ \gcd(k,N)=1}|\]
and by $\Psi(N)$ the two dimensional analog, i.e.
\[\Psi(N):=|\Set{(k_1,k_2) \in \Z^2}{0 \leq k_1,k_2 < N,\ \gcd(k_1,k_2,N)=1}|.\]
From the theorem above, together with the
well known results
\[\sum_{d|N}\phi(d)=N\ \mathrm{and}\ \sum_{d|N}\Psi(d)=N^2,\]
we easily obtain the following.

\begin{corollary}
Let $n \in \Z$ and $N \in \N$. If $N=1$, then $L(n,N)=0$. If $N
\not= 1$, then
\[L(n,N)= \frac{n(n+1)}{12}(\Psi(N)-3\phi(N)) + \frac{2}{3}\epsilon(n,N),\]
where
\[ \epsilon(n,N):=\left\{
\begin{array}{ll}
1 & \mbox{if $N=3$ and $n \equiv 1 \pmod{3}$}; \\
0 & \mbox{otherwise}.
\end{array}
\right. \]
\end{corollary}

\end{document}